\documentclass[%
reprint,
 amsmath,amssymb,
 aps,
]{revtex4-2}

\usepackage{graphicx}
\usepackage{dcolumn}
\usepackage{bm}
\usepackage{hyperref}
\usepackage{mathtools}
\usepackage{color}
\usepackage[dvipsnames]{xcolor}
\usepackage{amsmath,amsthm,float}
\usepackage{tikz}
\usepackage{tikz-cd}
\usetikzlibrary{matrix,arrows,decorations.pathmorphing,decorations.markings,shapes}
\usetikzlibrary{patterns} 

\newcommand{\R}[1]{\mathbb{R}^{#1}}

\newcommand{\Z}[1]{\mathbb{Z}^{#1}}

\newcommand{\bA}{\mathbf{A}}
\newcommand{\bb}{\mathbf{b}}
	
\newcommand{\bd}{\mathbf{d}}
\newcommand{\be}{\mathbf{e}}
\newcommand{\bE}{\mathbf{E}}
\newcommand{\bef}{\mathbf{f}}
\newcommand{\bF}{\mathbf{F}}

\newcommand{\bI}{\mathbf{I}}
\newcommand{\bJ}{\mathbf{J}}

\newcommand{\bv}{\mathbf{v}}
\newcommand{\bx}{\mathbf{x}}

\newcommand{\bphi}{\boldsymbol{\phi}}
\newcommand{\bPhi}{\boldsymbol{\Phi}}

\newcommand{\bzero}{\boldsymbol{0}}
\newcommand{\bcR}{\boldsymbol{\mathcal{R}}}

\newcommand{\balpha}{\boldsymbol{\alpha}}

\newcommand{\de}{\mathrm d}

\newcommand{\eps}{\varepsilon}

\newcommand{\average}{{\mathchoice {\kern1ex\vcenter{\hrule
height.4pt width 8pt depth0pt}
\kern-11pt} {\kern1ex\vcenter{\hrule height.4pt width 4.3pt
depth0pt} \kern-7pt} {} {} }}

\newtheorem{defin}{Definition}

\newtheorem{theorem}[defin]{Theorem}

\mathchardef\emptyset="001F

\begin{document}

\preprint{APS/123-QED}

\title{Controlling non-controllable scallops}

\author{Marta Zoppello and Marco Morandotti}%
\email{$\{$marta.zoppello,marco.morandotti$\}$@polito.it}
\affiliation{Dipartimento di Scienze Matematiche ``G.~L.~Lagrange'', Politecnico di Torino}%
\author{Hermes Bloomfield-Gad\^{e}lha}
\email{hermes.gadelha@bristol.ac.uk}
\affiliation{Department of Engineering Mathematics and Bristol Robotics Laboratory, University of Bristol.}

\date{\today}

\begin{abstract}
Any swimmer embedded on a inertialess fluid must perform a non-reciprocal motion to swim forward. The archetypal demonstration of this unique motion-constraint was introduced by Purcell with the so-called “scallop theorem”. Scallop here is a minimal mathematical model of a swimmer composed by two arms connected via a hinge whose  periodic motion (of opening and closing its arms) is not sufficient to achieve net displacement. Any source of incongruence on the motion or in the forces/torques experienced by such time-reversible scallop will break the symmetry imposed by the Stokes linearity and lead to subsequent propulsion of the scallop. However, little is known about the controllability of time-reversible scalloping systems. Here, we consider two individually non-controllable scallops swimming together. Under a suitable geometric assumption on the configuration of the system, it is proved that non-zero net displacement can be achieved as a consequence of their hydrodynamic interaction. A detailed analysis of the control system of equations is carried out analytically by means of geometric control theory. We obtain an analytic expression for the the displacement after a prescribed sequence of controls in function of the phase difference of the two scallops. Numerical validation of the theoretical results is presented with model predictions in further agreement with the literature.
 
\end{abstract}

\maketitle


\section{Introduction} \label{sec:intro}
Inertialess hydrodynamics is notorious for its time-reversibility constraint \cite{lighthill1975mathematical,Purcell1977}. Linearity of Stokes equation imprisons any swimmer moving in a time-reversible manner in perpetuity, a consequence of the “scallop theorem”. Purcell \cite{Purcell1977} introduced this notion by considering the simplest mathematical abstraction of a time-reversible swimmer that goes nowhere in a viscous fluid, the “scallop”, composed by two rigid arms that open and close relative to a hinge point, in an analogy to the opening movement scallops in the sea. Most interestingly, perhaps it is not accidental the endemic number of naturally occurring non-reciprocal swimmers in the microscopic world \cite{velho2021bank,lisicki2019swimming,gaffney2011mammalian,rossi2017kinematics}. From eukaryotic to prokaryotic microscopic organisms, the non-reciprocity in their motion in time is at the heart of their biological function, with elasticity most commonly explored by eukaryotic cilia and flagella, such for spermatozoa swimming \cite{gaffney2011mammalian} or ciliary beating \cite{marumo2021three}. 

With recent advances in manufacturing, prototyping and 3D printing, robotic swimming is an emerging field able to realise mathematical models that only lived in the abstract sense \cite{milana2020metachronal,sareh2013swimming,gu2020magnetic,qiu2014artificial,qiu2014swimming,dillinger2021starfish}. As such, recent research has been increasingly interested in the connection between swimming and control theory in viscous fluid environments \cite{lauga2008no,ADSL2008,ADSGZ2013,maggistro2019optimal}. Using Purcell’s archetypical swimmer, several studies focused their efforts in finding ways to break with the scallop theorem and enable Purcell’s scallop to move forward \cite{lauga2011life,ADSGZ2013,GMZ2015}. Breaking the time-reversibility constraint of a two-rigid arms scallop requires other sources of asymmetries in the system, such as higher degrees of freedom, as seen for Purcell’s three-link swimmer \cite{Purcell1977,BKS2003,ADSGZ2013,BBGMR2017}, non-linearities arising from complex fluid rheology \cite{lauga2009life,qiu2014swimming,pak2010pumping,yu2007experiments}, inertial effects \cite{bruot2021direct,lauga2011life,hubert2021scallop}, non-symmetric hydrodynamic interactions \cite{lauga2008no,lauga2011life,takagi2015swimming,reinmuller2013se}, stochasticity \cite{lauga2011enhanced}, among other interactions \cite{teboul2019breakdown}. Indeed, even though one scallop cannot swim forward over one period, two scallops swimming together can. This is due to incongruences on the forces/torques experienced by the swimmers via non-local hydrodynamic interactions \cite{MKL2016,TCL2021}, despite their individual reciprocal motion. This in turns allows the forward motion of multiple scallops swimming collectively \cite{lauga2008no}. Against this background, and despite of the large body of mathematical investigations focusing on how to break time-reversibility constraint of Purcell’s scallop, little is known about the controllability of this system under the influence of non-local hydrodynamic interactions. In the realms of Geometric Control Theory, controllability is defined as the existence of \emph{control functions} that are able to steer the system from a given initial configuration to a given final one \cite{AS2004,Coron2007}. More specifically, in the context of micro-swimmers, it is the ability of finding a suitable swimming stroke that generates a non-zero net displacement. 

In this work, we explore the controllability of two inherently non-controllable units by exploiting the non-local nature of their mutual hydrodynamic interaction. To allow analytical progress, we study the synchronized motion of two Purcell scallops within a minimal hydrodynamic interacting model between nearby discretized elements using control theory. We mathematically prove that the scallops’ shape shift determines the change in position and orientation of the system: thanks to the Invertibility Theorem~\ref{invertible}, we write the equations of motion as a control system; by tools from Geometric Control Theory, we obtain an analytic expression of the displacement after a prescribed sequence of controls, see formula \eqref{113}. Finally, we optimize the displacement as a function of the phase difference of the controls acting on each individual scallop and recover well-known results in the literature, see Theorem~\ref{taylor}. The numerical experiments undertaken further validate the analytical results and the theoretical predictions. Our novel result sets the basis for a deeper understanding of systems of mutually interacting micro-swimmers, whose choral motion is determined by the non-local hydrodynamic interactions generated by individual rate of  shape changes of each unit. A different approach using applied forces and torques as controls is proposed in \cite{walker2021control}; see also \cite{AY2008}, where the motion of two interacting dumb-bell swimmers is studied.


\section{The model}
We consider a system composed by two scallops of length $2L$ swimming with a distance $h$ from each other, as illustrated in Fig. \ref{2_link_sys}. We build our model on previous investigation by Man et al.~\cite{MKL2016} for two filaments coupled hydrodynamically at low Reynolds number, with setup adapted to our case. In particular, denoting by $a>0$ the thickness of the scallops and by $\sigma_i$ ($i=1,2$) is the opening angle of the $i$-th scallop, we assume that,
\begin{equation}\label{assumption}
a<<h<<L\quad\text{and}\quad \sigma_i\approx\pi.
\end{equation}
These assumptions are needed in order to use the approximation proposed in \cite{MKL2016} to compute the interaction forces and will allow us to compare, in Section~\ref{NV}, our numerical results with those presented in Ref.~\cite{MKL2016}. We notice that the second assumption in \eqref{assumption} means that the scallops are almost completely open, so that they can be considered approximately aligned. 

The system we consider is shown in Fig.~\ref{2_link_sys} and is described using the following coordinates: For $i=1,2$, let $t\mapsto \bx_i(t)\coloneqq(x_i(t),y_i(t))^\top$ denote the time-dependent position (with respect to the laboratory fixed frame of reference) of the hinge of $i$-th scallop, let $t\mapsto\theta_i(t)$ denote the time-dependent angle that the upper link of the $i$-th scallop forms with the positive $x$-axis, and let 
$$t\mapsto \be^{(j)}_i(t)\coloneqq 
\begin{pmatrix}\cos\big(\theta_i(t)+(j-1)\sigma_i(t)\big) \\[1mm] \sin\big(\theta_i(t)+(j-1)\sigma_i(t)\big)
\end{pmatrix}$$ 
denote the directions of the $j$-th link of the $i$-th scallop (for $j=1,2$), where $t\mapsto \sigma_i(t)$ is the opening angle of the $i$-th scallop. With this choice of coordinates, the triple $(x_i,y_i,\theta_i)$ is the set of coordinates describing the \emph{position} of the $i$-th scallop in the plane and $\sigma_i$ is its \emph{shape}. Finally, we let $L>0$ be the length of each link, so that the time-dependent position of the generic point $\bx^{(j)}_i(s,t)\in\R{2}$ on the $j$-th link of the $i$-th scallop at distance $s\in[0,L]$ from the hinge is given by 
\begin{equation}\label{100}
\bx^{(j)}_i(s,t)=\bx_i(t)+s\be^{(j)}_i(t)
\end{equation}
and its velocity is given by
\begin{equation*}
\begin{split}
\bv^{(j)}_i(s,t)\coloneqq &\,  \dot\bx_i^{(j)}(s,t)=\dot\bx_i(t)+s\dot\be_{i}^{(j)}(t) \\
=& \begin{pmatrix}\dot x_i(t)\\[1mm] \dot y_i(t)\end{pmatrix}+s\big(\be^{(j)}_i(t)\big)^\bot(\dot\theta_i(t)+(j-1)\dot\sigma_i(t)).
\end{split}
\end{equation*}
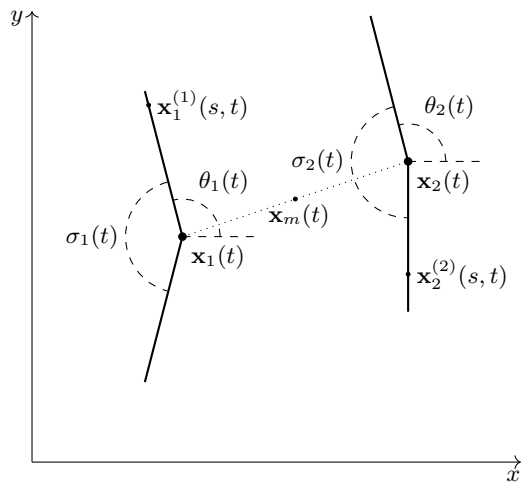
\begin{figure}[!h]
\begin{center}
\begin{tikzpicture}
\draw[->, color=black] (0,0) -- (6.5,0);
\node[text=black,below] at (6.4,0) {$x$};
\draw[->, color=black] (0,0) -- (0,6);
\node[text=black,left] at (0,5.9) {$y$};
\draw[color=black,fill] (2,3) circle (.05) (2,3) node [text=black,below right] {$\bx_1(t)$};
\draw[color=black,thick] (2,3) -- (1.5,4.936);
\draw[color=black,thick] (2,3) -- (1.5,1.064);
\draw[color=black,dashed] (2,3) -- (3,3);
\draw[color=black,dashed] (2.5,3) arc (0:105:.5) (3,4) node [text=black,below left] {$\theta_1(t)$}; 
\draw[color=black,dashed] (1.8,3.725) arc (105:255:.75) (.8,3) node [text=black] {$\sigma_1(t)$}; 
\draw[color=black,fill] (5,4) circle (.05) (5,4) node [text=black,below right] {$\bx_2(t)$};
\draw[color=black,fill] (1.55,4.75) circle (.025) (1.55,4.75) node [text=black,right] {$\bx_1^{(1)}(s,t)$};
\draw[color=black,dashed] (5,4) -- (6,4);
\draw[color=black,thick] (5,4) -- (5,2);
\draw[color=black,thick] (5,4) -- (4.5,5.936);
\draw[color=black,dashed] (5.5,4) arc (0:105:.5) (6,5) node [text=black,below left] {$\theta_2(t)$}; 
\draw[color=black,dashed] (4.8,4.725) arc (105:270:.75) (3.8,4) node [text=black] {$\sigma_2(t)$}; 
\draw[color=black,fill] (5,2.5) circle (.025) (5,2.5) node [text=black,right] {$\bx_2^{(2)}(s,t)$};
\draw[color=black,dotted] (2,3) -- (5,4);
\draw[color=black,fill] (3.5,3.5) circle (.025) (3.55,3.55) node [text=black,below] {$\bx_m(t)$};
\end{tikzpicture}
\caption{Schematic representation of two scallops swimming together relative to the laboratory fixed frame of reference.}
\label{2_link_sys}
\end{center}
\end{figure}

Denoting by $(s,t)\mapsto \bef^{(j)}_i(s,t)$ the density of hydrodynamic force acting on $\bx^{(j)}_i(s,t)$, Resistive Force Theory \cite{GH1955} states that it is proportional to the local velocities $\bv^{(j)}_i(s,t)$ of the swimmer relative the background fluid; taking also the hydrodynamics interaction between the two scallops into account, as studied in \cite{MKL2016}, we have, for $i,j=1,2$,
\begin{equation}\label{force_density}
\bef^{(j)}_i(s,t)=\frac{-1}{\Lambda(s,t)}\bJ^{(j)}_i(t) \bv^{(j)}_i(s,t)+\frac{\lambda(s,t)}{\Lambda(s,t)}\bJ^{(j)}_{\neg i}(t) \bv^{(j)}_{\neg i}(s,t)
\end{equation} 
where $\neg i\coloneqq 3-i$ is a concise form to write the value not taken by index $i$, $\lambda(s,t)\coloneqq \displaystyle \frac{\ln \big(h(s,t)/L\big)}{\ln (a/L)}\in(0,1)$ (see Section~\ref{rem_lambda} below for more comments), and $\Lambda(s,t)\coloneqq 1-\lambda^2(s,t)$. In formula \eqref{force_density}, $t\mapsto \bJ^{(j)}_i(t)$ is the Resistive Force Theory operator relative to the $j$-th link of the $i$-th swimmer and it is defined by
\begin{equation}\label{103}
\bJ^{(j)}_i(t)=
C_\perp\bI+(C_\parallel-C_\perp)\bE_i^{(j)}(t), 
\end{equation}
with $\bE_i^{(j)}(t)\coloneqq \be^{(j)}_i(t)\otimes\be^{(j)}_i(t)$. 
Notice that, by the first assumption in \eqref{assumption}, it is reasonable to assume that the $j$-th link of each scallop interacts only with the $j$-th link of the other scallop, neglecting the interaction with the $\neg j$-th link. This is due to the fact that, since $h<<L$, the scallops can only vary their shape minimally to avoid overlapping.
In \eqref{103}, the constants $C_\perp$ and $C_\parallel$ are the drag coefficients in the perpendicular and parallel directions, respectively, to the links and are measured in $\mathrm{Ns}/\mu\mathrm{m}^2$~\footnote{We use SI units throughout the paper}. In general, for micro-swimmers, a physically meaningful assumption is to take $C_\perp\approx2C_\parallel$. 
By the second assumption in \eqref{assumption}, 
it is not restrictive to suppose that $h(s,t)$ does not depend on spatial variable $s$ and that it undergoes very little variations in time, so that it can be considered constant, that is $h(s,t)=h$; thus, we can also consider $\lambda(s,t)=\lambda$ and $\Lambda(s,t)=\Lambda=1-\lambda^2$.

The time-dependent force acting on the $j$-th link of the $i$-th scallop is obtained by integration of \eqref{force_density} over $s\in[0,L]$, namely
\begin{equation}\label{101}
\bF_i^{(j)}(t)=\int_0^L \bef_i^{(j)}(s,t)\,\de s,
\end{equation}
so that the overall time-dependent force $t\mapsto \bF_i(t)$ on the $i$-th scallop is given by $\bF_i(t)=\bF_i^{(1)}(t)+\bF_i^{(2)}(t)$, that is,
\begin{equation}\label{102}
\bF_i(t)=
\int_0^L \big(\bef_i^{(1)}(s,t)+\bef_i^{(2)}(s,t)\big)\,\de s.
\end{equation}
From \eqref{force_density} and \eqref{103} we have (dropping the dependences on~$s$ and~$t$)
\begin{equation*} 
\begin{aligned}
&\bef^{(1)}_i+\bef^{(2)}_i 
= -\frac{1}{\Lambda} \big[2C_\perp\bI+(C_\parallel-C_\perp) \big(\bE_i^{(1)}+\bE_i^{(2)}\big)\big]\dot\bx_i \\
&\, -\frac1\Lambda sC_\perp\big(\be_i^{(1)}+\be_i^{(2)}\big)^\perp\dot\theta_i-\frac1\Lambda sC_\perp \big(\be_i^{(2)}\big)^\perp \dot\sigma_i \\
&\, + \frac{\lambda}{\Lambda}\big[2C_\perp\bI+(C_\parallel-C_\perp) \big(\bE_{\neg i}^{(1)}+\bE_{\neg i}^{(2)}\big)\big]\dot\bx_{\neg i} \\
&\, +\frac{\lambda}{\Lambda} sC_\perp\big(\be_{\neg i}^{(1)}+\be_{\neg i}^{(2)}\big)^\perp\dot\theta_{\neg i}+\frac{\lambda}{\Lambda} sC_\perp \big(\be_{\neg i}^{(2)}\big)^\perp \dot\sigma_{\neg i}\,. 
\end{aligned}
\end{equation*}
Therefore formula \eqref{102} reads
\begin{equation}\label{104}
\begin{split}
&\bF_i(t)= -\frac1\Lambda \big(\bA_i(t)\dot\bx_i(t)+ \bb_i(t)\dot\theta_i(t)+\balpha_i(t)\dot\sigma_i(t)\big) \\
& +\frac{\lambda}{\Lambda} \big(\bA_{\neg i}(t)\dot\bx_{\neg i}(t) + \bb_{\neg i}(t)\dot\theta_{\neg i}(t)+\balpha_{\neg i}(t)\dot\sigma_{\neg i}(t)\big),
\end{split}
\end{equation}
where we have defined 
\begin{equation}\label{105}
\begin{split}
\bA_i(t)\coloneqq &\, L\big[2C_\perp\bI+(C_\parallel-C_\perp) \big(\bE_i^{(1)}(t)+\bE_i^{(2)}(t)\big)\big], \\
\bb_i(t)\coloneqq &\, \frac{L^2}2 C_\perp\big(\be_i^{(1)}(t)+\be_i^{(2)}(t)\big)^\perp, \\
\balpha_i(t)\coloneqq &\, \frac{L^2}2 C_\perp \big(\be_i^{(2)}(t)\big)^\perp.
\end{split}
\end{equation}
To complete the set of equation of motion, we also need to compute the torque acting on each swimmer. To this end, we start from the torque density with respect to the $i$-th hinge, which reads, using \eqref{100}, \eqref{force_density}, and \eqref{103},
\begin{equation*}\label{torque_density}
\begin{split}
&\tau^{(j)}_i(s,t)=  \big(\bx_i^{(j)}(s,t)-\bx_i(t)\big)\times\bef_i^{(j)}(s,t) \\ 
&\, = -\frac1{\Lambda} C_\perp s \big(\be_i^{(j)}\big)^\perp \!\! \cdot\dot\bx_i-\frac1\Lambda C_\perp s^2 \dot\theta_i+\frac\lambda\Lambda C_\perp s \big(\be_i^{(j)}\big)^\perp \!\! \cdot \dot\bx_{\neg i} \\
&\, +\frac\lambda\Lambda (C_\parallel-C_\perp)s\Big(\big(\be_i^{(j)}\big)^\perp\! \cdot \be_{\neg i}^{(j)}\Big)\big(\be_{\neg i}^{(j)} \!\cdot \dot \bx_{\neg i}\big) \\
&\, +\frac\lambda\Lambda C_\perp s^2  \big(\be_i^{(j)} \cdot \be_{\neg i}^{(j)}\big)\dot\theta_{\neg i} -\frac{j-1}\Lambda C_\perp s^2 \dot\sigma_{i} \\
&\, +(j-1)\frac\lambda\Lambda C_\perp s^2  \big(\be_i^{(j)}\cdot \be_{\neg i}^{(j)}\big)\dot\sigma_{\neg i}\,.
\end{split}
\end{equation*}

Analogously to what we did in \eqref{102}, the overall time-dependent torque $t\mapsto \tau_i(t)$ acting on the $i$-th scallop is given by $T_i(t)=T_i^{(1)}(t)+T_i^{(2)}(t)$, that is,
\begin{equation}\label{102t}
T_i(t)= \int_0^L \! \big(\tau_i^{(1)}(s,t)+\tau_i^{(2)}(s,t)\big)\,\de s,
\end{equation}
which has the expression
\begin{equation}\label{103t}
\begin{split}
\!\!\! T_i(t)=&\, -\frac1\Lambda \bb_i(t) \cdot \dot\bx_i(t)-\frac1\Lambda \omega \dot\theta_i(t)+\frac\lambda\Lambda \bd_i(t) \cdot \dot\bx_{\neg i}(t)\\
&\,+\frac\lambda\Lambda \varpi(t)\dot\theta_{\neg i}(t) -\frac1\Lambda \frac{\omega}2 \dot\sigma_i(t)+\frac\lambda\Lambda \beta(t) \dot\sigma_{\neg i}(t),
\end{split}
\end{equation}
where $\bb_i(t)$ is as in \eqref{105} and 
\begin{equation}\label{104t}
\begin{split}
\!\!\! \bd_i(t)\coloneqq&\, \frac{L^2(C_\parallel-C_\perp)}2\Big[\!\Big(\big(\be_i^{(1)}(t)\big)^\perp\! \cdot \be_{\neg i}^{(1)}(t)\Big)\be_{\neg i}^{(1)}(t) \\
&+\! \Big(\big(\be_i^{(2)}(t)\big)^\perp\! \cdot \be_{\neg i}^{(2)}(t)\Big)\be_{\neg i}^{(2)}(t)\Big] + \bb_i(t), \\
\omega \coloneqq&\, \frac{2L^3C_\perp}3, \\
\!\! \varpi(t)\coloneqq&\, \frac{L^3C_\perp}3 \big[\be_i^{(1)}(t)\cdot\be_{\neg i}^{(1)}(t)+ \be_i^{(2)}(t)\cdot\be_{\neg i}^{(2)}(t)\big], \\ 
\beta(t)\coloneqq&\, \frac{L^3C_\perp}3 \be_i^{(2)}(t)\cdot\be_{\neg i}^{(2)}(t), 
\end{split}
\end{equation}
where we notice that both $\varpi(t)$ and $\beta(t)$ do not depend on $i$ since they are symmetric in the exchange $i\mapsto\neg i$.
Formulae \eqref{104}, \eqref{105}, \eqref{103t}, and \eqref{104t} can be gathered together in the following expression
\begin{widetext}
\begin{equation}\label{106}
\begin{split}
-\Lambda&\begin{pmatrix}
\bF_1(t) \\
T_1(t) \\
\bF_2(t) \\
T_2(t)
\end{pmatrix}= 
\left(\begin{array}{cc|cc}
\bA_1(t) & \bb_1(t) & -\lambda\bA_2(t) & -\lambda\bb_2(t) \\
\bb_1^\top(t) & \omega & -\lambda\bd_1^\top(t) & -\lambda\varpi(t) \\
\hline 
-\lambda\bA_1(t) & -\lambda\bb_1(t) & \bA_2(t) & \bb_2(t) \\
-\lambda\bd_2^\top(t) & -\lambda\varpi(t) & \bb_2^\top(t) & \omega
\end{array}\right)\begin{pmatrix}
\dot\bx_1(t) \\
\dot\theta_1(t) \\
\dot\bx_2(t) \\
\dot\theta_2(t)
\end{pmatrix} 
+
\begin{pmatrix}
\balpha_1(t) \\
\omega/2 \\
-\lambda\balpha_1(t) \\
-\lambda\beta(t)
\end{pmatrix}\dot\sigma_1(t)
+
\begin{pmatrix}
-\lambda\balpha_2(t) \\
-\lambda\beta(t) \\
\balpha_2(t) \\
\omega/2
\end{pmatrix}\dot\sigma_2(t) \\
&\,\eqqcolon  
\bcR(t;\lambda)\!\begin{pmatrix}
\dot\bx_1(t) \\
\dot\theta_1(t) \\
\dot\bx_2(t) \\
\dot\theta_2(t)
\end{pmatrix} +
\bphi_1(t;\lambda)\dot\sigma_1(t)+
\bphi_2(t;\lambda)\dot\sigma_2(t) 
=
\bcR(t;\lambda)\!\begin{pmatrix}
\dot\bx_1(t) \\
\dot\theta_1(t) \\
\dot\bx_2(t) \\
\dot\theta_2(t)
\end{pmatrix} +
\big(\bphi_1(t;\lambda)|\bphi_2(t;\lambda)\big)\!\begin{pmatrix}
\dot\sigma_1(t) \\
\dot\sigma_2(t)
\end{pmatrix} \\
&\,= 
\bcR(t;\lambda)\begin{pmatrix}
\dot\bx_1(t) \\
\dot\theta_1(t) \\
\dot\bx_2(t) \\
\dot\theta_2(t)
\end{pmatrix} +
\bPhi(t;\lambda)\begin{pmatrix}
\dot\sigma_1(t) \\
\dot\sigma_2(t)
\end{pmatrix}.
\end{split}
\end{equation}
\end{widetext}
\section{The equations of motion}
The equations of motion, which are obtained by imposing the total force and torque balance, that is, that the left-hand side in \eqref{106} is equal to zero, read
\begin{equation}\label{107}
\bcR(t;\lambda)\begin{pmatrix}
\dot\bx_1(t) \\
\dot\theta_1(t) \\
\dot\bx_2(t) \\
\dot\theta_2(t)
\end{pmatrix} + \bPhi(t;\lambda)\begin{pmatrix}
\dot\sigma_1(t) \\
\dot\sigma_2(t)
\end{pmatrix}=\bzero.
\end{equation}
The first step to solve \eqref{107} is to investigate the invertibility of the matrix 
\begin{equation}\label{RRR}
\bcR(t;\lambda)=\begin{pmatrix}
\bcR_{11}(t) & -\lambda\bcR_{12}(t) \\
-\lambda \bcR_{21}(t) & \bcR_{22}(t)
\end{pmatrix},
\end{equation}
where we observe that the blocks $\bcR_{11}$ and $\bcR_{22}$ (not depending on $\lambda$) are the grand resistance matrices \cite{HB1965} of the individual scallops and so they are both (positive) definite and symmetric, and therefore invertible. 
Therefore, we have following theorem.
\begin{theorem}[Invertibility]\label{invertible}
There exists $\lambda_0\in(0,1)$ such that the matrix $\bcR(t;\lambda)$ defined in \eqref{RRR} is invertible for every $\lambda\in[0,\lambda_0)$ and for every $t\in[0,+\infty)$.
\end{theorem}
\begin{proof}
Since $\bcR_{11}(t)$ and $\bcR_{22}(t)$ are invertible, we have 
\begin{equation}\label{det0}
\begin{split}
\det \big(\bcR(t;0)\big)=&\det\begin{pmatrix}
\bcR_{11}(t) & \bzero_{3\times3} \\
\bzero_{3\times3} & \bcR_{22}(t)
\end{pmatrix}\\
=&\det\big(\bcR_{11}(t)\big)\det\big(\bcR_{22}(t)\big)\neq0
\end{split}
\end{equation}
for every $t\in[0,+\infty)$. 
Recalling that $\det\big(\bcR(t;\lambda)\big)=\det\big(\bcR_{11}(t)\big)\det\big(\bcR_{22}(t)-\lambda^2\bcR_{12}(t)\bcR_{11}^{-1}(t)\bcR_{21}(t)\big)$ (see \cite[Chapter 7.7]{HJ2013})
and noticing that the expression is continuous in $\lambda$, \eqref{det0} implies that there exists a value $\lambda_0\in(0,1)$ such that $\det \big(\bcR(t;\lambda)\big)\neq0$ for $\lambda\in[0,\lambda_0)$, so that $\bcR(t;\lambda)$ is invertible for $(t,\lambda)\in[0,+\infty)\times[0,\lambda_0)$.
\end{proof}

We propose two approaches to study the equations of motion~\eqref{107}. 
The first is by using standard results from ODE theory; from this point of view, isolating the contribution of the shape parameters $\bPhi(\dot\sigma_1,\dot\sigma_2)^\top$ embodies the ability of swimming by shape deformation, as it will be clear from \eqref{108} below.
The latter is by relying on control theory, where the velocities $(\dot\sigma_1,\dot\sigma_2)^\top$ can be considered as the controls of the system. Ideally, those are actuators that can be prescribed to steer the system. In particular, it will be possible to quantify the motion in terms of given controls.

\subsection{Approach by ODE theory}
By solving \eqref{107} for the translational and rotational velocities, by Theorem~\ref{invertible} we obtain 
\begin{equation}\label{108}
\begin{pmatrix}
\dot\bx_1(t) \\
\dot\theta_1(t) \\
\dot\bx_2(t) \\
\dot\theta_2(t)
\end{pmatrix}=-\bcR(t;\lambda)^{-1}\bPhi(t;\lambda)\begin{pmatrix}
\dot\sigma_1(t) \\
\dot\sigma_2(t)
\end{pmatrix},
\end{equation}
expressing the fact that the position and orientation of the scallops are determined by their shape change. Invoking standard results on ordinary differential equations \cite{Hale1980} (see, \emph{e.g.}, \cite{ADSL2008,ADSGZ2013} for the three-sphere swimmer and the $N$-link swimmer; see also \cite[Theorem~6.4]{DMDSM2011} for the abstract setting and \cite[Theorem~3.3]{DMDSM2013} for the case of a planar one-dimensional swimmer), once the deformation $t\mapsto(\sigma_1(t),\sigma_2(t))$ is prescribed and satisfies certain regularity conditions, for a given initial datum $(\bx_{1}^\circ,\theta_1^\circ,\bx_2^\circ,\theta_2^\circ)$, the initial-value problem for system \eqref{108} admits a unique solution, which also continuously depends on the initial data.

\subsection{Approach by control theory}
Equation \eqref{107} can also be interpreted as a control system by setting $(\dot\sigma_1,\dot\sigma_2)^\top=(u_1,u_2)^\top$, for suitable control maps $t\mapsto(u_1(t),u_2(t))$. With this notation, \eqref{107} becomes
\begin{equation}\label{109}
\!\!\!\! \!\!
\left(\!\begin{array}{c|c}
\bcR(t;\lambda) & \bzero_{6\times2} \\
\hline
\bzero_{2\times6} & \bI_{2\times2}
\end{array}\right) \!\!
\left(\begin{array}{c}
\dot\bx_1(t) \\
\dot\theta_1(t) \\
\dot\bx_2(t) \\
\dot\theta_2(t) \\
\hline
\dot\sigma_1(t) \\
\dot\sigma_2(t)
\end{array}\!\right) \!\!= \!\! \left(\!\begin{array}{c}
-\bPhi(t;\lambda) \\
\hline
\bI_{2\times2}
\end{array}\!\right) \!\! \begin{pmatrix}
u_1(t) \\
u_2(t)
\end{pmatrix}.
\end{equation}
Invoking the Invertibility Theorem~\ref{invertible} again, 
we have
\begin{equation}\label{110}
\begin{split}
\!\!\!\!
\left(\begin{array}{c}
\dot\bx_1(t) \\
\dot\theta_1(t) \\
\dot\bx_2(t) \\
\dot\theta_2(t) \\
\hline
\dot\sigma_1(t) \\
\dot\sigma_2(t)
\end{array}\right)=&\, \left(\begin{array}{c}
-\bcR(t;\lambda)^{-1}\bPhi(t;\lambda) \\
\hline
\bI_{2\times2}
\end{array}\right)\begin{pmatrix}
u_1(t) \\
u_2(t)
\end{pmatrix} \\
=&\, \left(\begin{array}{c}
-\bcR(t;\lambda)^{-1}\bphi_1(t;\lambda) \\
\hline
1 \\
0
\end{array}\right)u_1(t) \\
&\,+\left(\begin{array}{c}
-\bcR(t;\lambda)^{-1}\bphi_2(t;\lambda) \\
\hline
0 \\
1
\end{array}\right)u_2(t) \\[2mm]
\eqqcolon&\, \bv_1(\theta_1(t),\theta_2(t),\sigma_1(t),\sigma_2(t);\lambda)u_1(t) \\
&\, +\bv_2(\theta_1(t),\theta_2(t),\sigma_1(t),\sigma_2(t);\lambda)u_2(t),
\end{split}
\end{equation}
which is a drift-less affine control system. System \eqref{110} can be studied from the point of view of Geometric Control Theory \cite{AS2004,Coron2007}. We pursue this analysis in Section~\ref{BoS}, where we show, by computing the Lie brackets $[\bv_1,\bv_2]$, that the system undergoes a non-zero net displacement. 
In studying some special deformations $t\mapsto(u_1(t),u_2(t))^\top$, we also show numerically in Section~\ref{NV} that there exists a choice of $u_1$ and $u_2$ such that the displacement is maximized.

\section{Breaking of symmetry}\label{BoS}
Purcell's celebrated Scallop Theorem \cite{Purcell1977} states that if the two $2$-links were considered individually, they would not be able to achieve a non-zero net displacement as a consequence of a reciprocal motion: just by opening and closing the angle $\sigma_i$ in a periodic way neither $2$-link can advance. From the geometrical control theory viewpoint, this is due to the fact that one vector field alone cannot generate enough degrees of freedom for the $2$-links to advance.
For this reason, a system of two $2$-links which beat synchronously, i.e., with $t\mapsto\sigma_1(t)=\sigma_2(t)$, still cannot advance since there is only one shape variable.

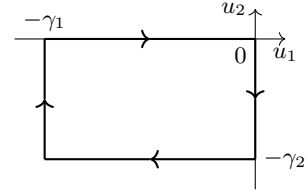
\begin{figure}[!h]
\begin{center}
\begin{tikzpicture}[scale=.8]
\draw[->, color=black] (-4,0) -- (.5,0);
\draw[->, color=black] (0,-2.5) -- (0,.5);
\node[text=black,left] at (0,.5) {$u_2$};
\node[text=black,below] at (.5,0) {$u_1$};
\node[text=black,below left] at (0,0) {$0$};
\draw[->, color=black,thick] (0,0) -- (0,-1);
\draw[color=black,thick] (0,-1) -- (0,-2);
\draw[->,color=black,thick] (0,-2) -- (-1.75,-2);
\draw[color=black,thick] (-1.75,-2) -- (-3.5,-2);
\draw[->, color=black,thick] (-3.5,-2) -- (-3.5,-1);
\draw[color=black,thick] (-3.5,-1) -- (-3.5,0);
\draw[->,color=black,thick] (-3.5,0) -- (-1.75,0);
\draw[color=black,thick] (-1.75,0) -- (0,0);
\node[text=black,right] at (0,-2) {$-\gamma_2$};
\node[text=black,above] at (-3.5,0) {$-\gamma_1$};
\end{tikzpicture}
\caption{The graph of the control stroke in \eqref{112}.}
\label{mettere}
\end{center}
\end{figure}

We now propose a way to overcome the Scallop Theorem by considering two time-dependent maps $t\mapsto\sigma_1(t)$ and $t\mapsto\sigma_2(t)$ different from each another. 
To start with, recalling \eqref{assumption}, we consider an initial configuration which is a perturbation of the aligned one (given by $\theta_i=\theta\in[0,2\pi)$ and $\sigma_i=\pi$, for $i=1,2$, see Figure~\ref{2_link_sys}), namely
\begin{equation}\label{111}
\theta_i^\circ=\theta^\circ,\quad\sigma_i^\circ=\pi+\eps\cos((i-1)\phi),
\end{equation}
for $i=1,2$, $\eps>0$ a small parameter, and $\phi$ a phase. 
Next, we prescribe the following stroke in the time interval $[0,4\tau]$, for a small $\tau>0$ and for $\gamma_1,\gamma_2>0$,
\begin{equation}\label{112}
t\mapsto \begin{pmatrix}u_1(t)\\ u_2(t)\end{pmatrix}\coloneqq\begin{cases}
(0,-\gamma_2)^\top & \text{for $t\in[0,\tau)$,} \\
(-\gamma_1,0)^\top & \text{for $t\in[\tau,2\tau)$,} \\
(0,\gamma_2)^\top & \text{for $t\in[2\tau,3\tau)$,} \\
(\gamma_1,0)^\top & \text{for $t\in[3\tau,4\tau)$,} 
\end{cases}
\end{equation}
which corresponds to running clockwise along the boundary of the rectangle with horizontal side $\gamma_1$ and vertical side $\gamma_2$ which is located in the third quadrant in the $u_1u_2$-plane (see Figure~\ref{mettere}).
We now use \eqref{112} to compute the Lie brackets;
classical tools in geometric control theory \cite{Coron2007} yield that the solution to system \eqref{110} with initial conditions \eqref{111} is given by
\begin{equation}\label{113}
\begin{split}
&\!\!\!\! \left(\begin{array}{c}
\bx_1(4\tau) \\
\theta_1(4\tau) \\
\bx_2(4\tau) \\
\theta_2(4\tau) \\
\hline
\sigma_1(4\tau) \\
\sigma_2(4\tau)
\end{array}\right)=\left(\begin{array}{c}
\bx_1^\circ \\
\theta^\circ \\
\bx_2^\circ \\
\theta^\circ \\
\hline
\sigma_1^\circ \\
\sigma_2^\circ
\end{array}\right)-\gamma_1\gamma_2\tau^2 \bv_3^\circ(\phi;\theta^\circ) 
+o(\tau^2) \\
\!\!\!\! =&\! \left( \! \begin{array}{c}
\bx_1^\circ \\
\theta^\circ \\
\bx_2^\circ \\
\theta^\circ \\
\hline
\sigma_1^\circ \\
\sigma_2^\circ
\end{array} \! \right) \! -\gamma_1\gamma_2 \tau^2 \!  \left( \! \begin{array}{c}
\xi_1(\phi;\theta^\circ)\eps^2+ o(\eps^2) \\
\eta_1(\phi;\theta^\circ)\eps^2 +o(\eps^2) \\
\vartheta_1(\phi)\eps+o(\eps^2) \\
\xi_2(\phi;\theta^\circ)\eps^2+ o(\eps^2) \\
\eta_2(\phi;\theta^\circ)\eps^2 +o(\eps^2) \\
\vartheta_2(\phi)\eps+o(\eps^2) \\
\hline
0 \\
0
\end{array} \! \right) \! +o(\tau^2)
\end{split}
\end{equation}
where we have expanded the Lie bracket $\bv_3^\circ(\phi;\theta^\circ) 
\coloneqq [\bv_1(\cdot;\lambda),\bv_2(\cdot;\lambda)]|_{(\theta^\circ,\theta^\circ,\pi+\eps,\pi+\eps\cos\phi)}$ in powers of $\eps$ (up to second order) with
\begin{equation}\label{infinitesimal}
\begin{split}
\!\! \xi_1(\phi;\theta^\circ) \! =&\, \frac{L\lambda\cos\theta^\circ\sin^2\frac\phi2}{64C_\parallel C_\perp(1-\lambda^2)}\big((C_\perp^2(2+\lambda)-C_\parallel C_\perp) \\
&\,+\cos\phi(3C_\parallel^2-2C_\perp C_\parallel-C_\perp^2)\big) \\
\!\! \eta_1(\phi;\theta^\circ) \! =&\, \frac{L\lambda\sin\theta^\circ\sin^2\frac\phi2}{64C_\parallel C_\perp(1-\lambda^2)}\big((C_\perp^2(2+\lambda)-C_\parallel C_\perp) \\
&\, +\cos\phi(3C_\parallel^2-2C_\perp C_\parallel-C_\perp^2)\big) \\
\vartheta_1(\phi) \! =&\, \vartheta_2(\phi)=-\frac{\lambda}{16(1-\lambda)}\sin^2\frac\phi2 \\
\!\! \xi_2(\phi;\theta^\circ) \! =&\, \frac{L\lambda\cos\theta^\circ\sin^2\frac\phi2}{64C_\parallel C_\perp(1-\lambda^2)}\big((3C_\parallel^2-2C_\parallel C_\perp-C_\perp^2) \\
&\, +\cos\phi(C_\perp^2(2+\lambda)-C_\perp C_\parallel)\big) \\
\!\! \eta_2(\phi;\theta^\circ) \! =&\, \frac{L\lambda\sin\theta^\circ\sin^2\frac\phi2}{64C_\parallel C_\perp(1-\lambda^2)}\big((3C_\parallel^2-2C_\parallel C_\perp-C_\perp^2) \\
&\, +\cos\phi(C_\perp^2(2+\lambda)-C_\perp C_\parallel)\big).
\end{split}
\end{equation}

We notice that the bracket $\bv_3^\circ(\phi;\theta^\circ)$ 
in \eqref{113} is different from the zero vector, so there is a non-zero net displacement. Moreover, we notice that the initial shape has been restored, as expected from the periodicity of the motion (compare the last two components in \eqref{113} with \eqref{111}).

The two scallops rotate counter-clockwise by the same amount, which is of order $\eps$; for a clockwise rotation, it suffices to change the sign of either $\gamma_1$ or $\gamma_2$.
There is a net motion of order $\eps^2$ along both axes, which vanishes (up to order $o(\eps^2)$) according to the value of $\theta^0$: for instance, if $\theta^\circ=\pi/2$ the motion along the $x$-axis is negligible.
From \eqref{113} and \eqref{infinitesimal} we can estimate the global net displacement of the system by tracking the midpoint $t\mapsto\bx_m(t)=(x_m(t),y_m(t))$ of the line connecting the two hinges (see Figure~\ref{2_link_sys}). We have, up to $o(\tau^2)$,
\begin{equation}\label{midpoint}
\begin{split}
& \Delta\bx_m=\bx_m(4\tau)-\bx_m^\circ \\
&= \frac{\bx_1(4\tau)+\bx_2(4\tau)-\bx_1^\circ-\bx_2^\circ}2 \\
&= \frac{-\gamma_1\gamma_2\tau^2}2\begin{pmatrix}
(\xi_1(\phi;\theta^\circ)+\xi_2(\phi;\theta^\circ))\eps^2+o(\eps^2) \\
(\eta_1(\phi;\theta^\circ)+\eta_2(\phi;\theta^\circ))\eps^2+o(\eps^2) 
\end{pmatrix} \\
&= \frac{-\gamma_1\gamma_2\tau^2}2\begin{pmatrix}
C\eps^2\cos\theta^\circ\sin^2\phi+o(\eps^2) \\
C\eps^2\sin\theta^\circ\sin^2\phi+o(\eps^2) 
\end{pmatrix},
\end{split}
\end{equation}
where $C=C(L,\lambda,C_\parallel,C_\perp)$ is given by
\begin{equation}\label{C100}
C=
\frac{L\lambda \big(C_\perp^2(1+\lambda)-3C_\parallel C_\perp+3C_\parallel^2\big)}{128C_\parallel C_\perp(1-\lambda^2)}.
\end{equation}
Since $\lambda\in(0,1)$ and commonly for slender micro-swimmers one can take $C_\perp\approx 2C_\parallel$, which yields
\begin{equation}\label{C101}
\widetilde{C}\coloneqq C(L,\lambda,C_\parallel,2C_\parallel)=\frac{L\lambda(1+4\lambda)}{256(1-\lambda^2)}>0,
\end{equation}
it is easy to see that the constant $C$ in \eqref{C100} remains positive for values of $C_\perp$ sufficiently close to $2C_\parallel$.
Therefore, the net displacement of the midpoint $\bx_m$ is given, from \eqref{midpoint}, by
\begin{equation}\label{midpointnet}
\delta_m(\phi)\coloneqq|\Delta\bx_m|=
\frac{C\gamma_1\gamma_2\tau^2\eps^2\sin^2\phi}{2}(1+o(\eps^2)),
\end{equation}
which, in turn, can be maximized with respect to the phase $\phi$; it is easy to see that $\delta_m(\phi)$ is maximum for $\phi=\pi/2+k\pi$, for $k\in\Z{}$. 
Thus we have proved the following theorem, which recovers the classical result obtained by \cite{Taylor1951} (see also \cite{MKL2016}).
\begin{theorem}\label{taylor}
The maximal displacement for the pair of scallops is obtained for strokes that  have a phase difference of $\phi=\pi/2$. \qed
\end{theorem}

\section{Numerical validation}\label{NV}
\begin{figure}[H]
\centering
\includegraphics[scale=.275]{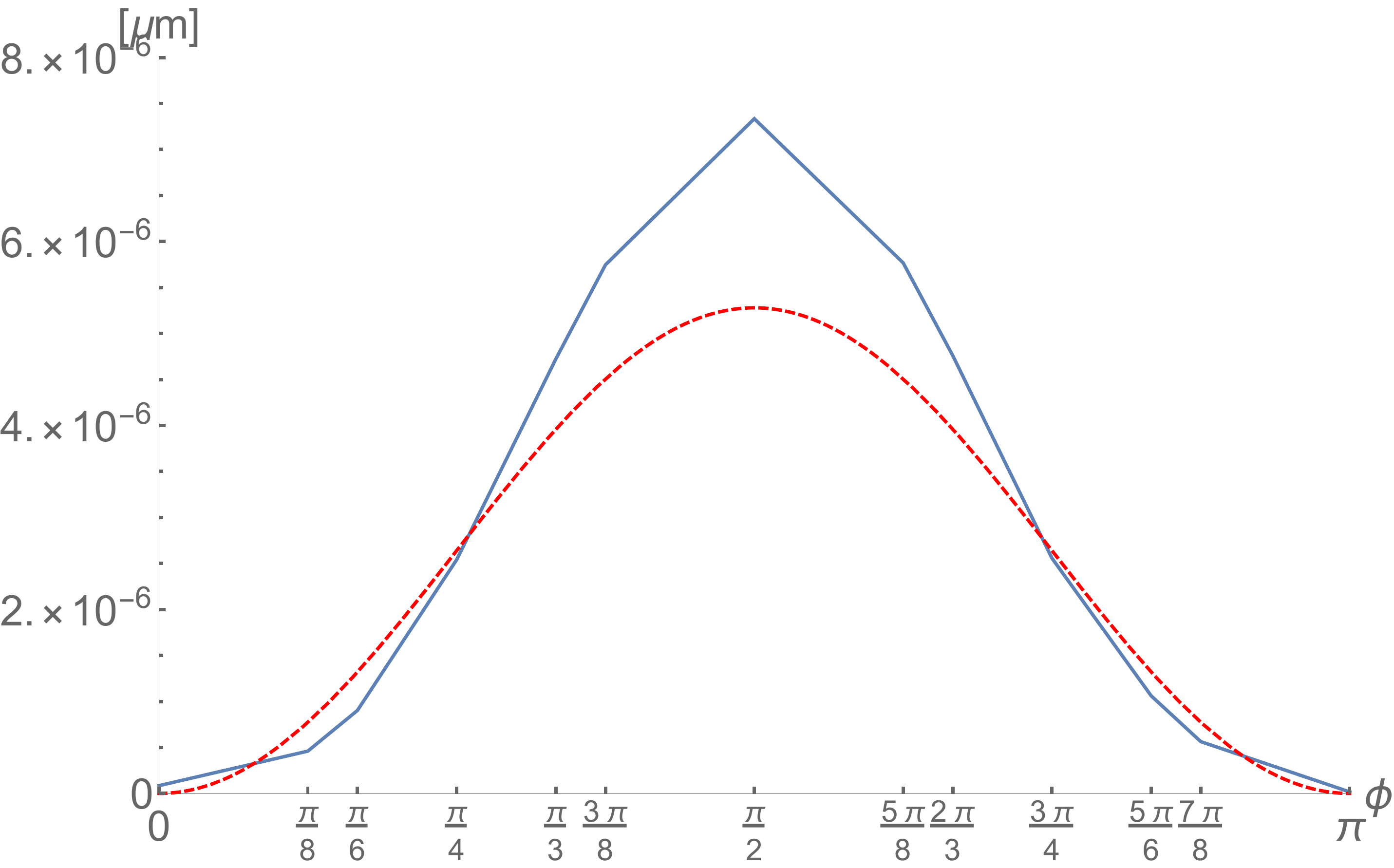}
\caption{Plot of $\delta_m^{(N)}=\big|\Delta\bx_m^{(N)}\big|$ as a function of $\phi$ (solid blue line) versus the theoretical $\delta_m(\phi)$ from \eqref{midpointnet} (dashed red line).}
\label{Delta_Xm}
\end{figure}

The second equality in \eqref{111} is the evaluation at $t=0$ of the map
\begin{equation}\label{111b}
\begin{split}
t\mapsto \sigma_i(t)=&\, \pi+\eps\cos\big(\omega t+(i-1)\phi\big) \\
=&\, \pi+\eps\cos\bigg(\frac{\pi t}{2\tau}+(i-1)\phi\bigg),
\end{split}
\end{equation}
which we can use to integrate numerically the equations of motion. The choice of the frequency $\omega=\pi/2\tau$ is made so that in the time interval $[0,4\tau]$ the angles $\sigma_i$ have returned to their initial value after spanning only one period. 
We have performed the numerical integration with the following set of parameters: $\eps=0.1$, $\omega=20$,
\begin{equation}\label{ni001}
\begin{split}
& L=10\,\mu\textrm{m},\quad h= 1\,\mu\textrm{m},\quad a=0.25\,\mu\textrm{m},\\ 
& C_\perp=2C_\parallel=2\,\textrm{Ns}/\mu\textrm{m}^2,
\end{split}
\end{equation}
for the following values of the phase $\phi$
\begin{equation*}
0,\, \frac\pi8,\, \frac\pi6,\, \frac\pi4,\, \frac\pi3,\, \frac{3\pi}8,\, \frac\pi2,\, \frac{5\pi}8,\, \frac{2\pi}3,\, \frac{3\pi}4,\, \frac{5\pi}6,\, \frac{7\pi}8,\, \pi.
\end{equation*}
Denoting by $\Delta\bx_m^{(N)}$ the displacement of the midpoint defined in \eqref{midpoint} evaluated numerically, we can plot the piecewise linear interpolation of the magnitude $\delta_m^{(N)}\coloneqq\big|\Delta\bx_m^{(N)}\big|$ as a function of $\phi$, obtaining the graph in Figure~\ref{Delta_Xm}, from which it is evident that $\delta_m^{(N)}$ is maximum for $\phi=\pi/2$, namely when the two filaments beat out of phase.

It is possible to compare $\delta_m=|\Delta\bx_m|$ and $\delta_m^{(N)}=\big|\Delta\bx_m^{(N)}\big|$, when $\phi=\pi/2$, upon choosing the appropriate values of~$\gamma_1$ and~$\gamma_2$ in \eqref{112}. Considering that the controls $t\mapsto(u_1(t),u_2(t))^\top$ are the derivatives of the angles $t\mapsto(\sigma_1(t),\sigma_2(t))^\top$, we can use \eqref{111b} (with $\phi=\pi/2$) to get
\begin{equation}\label{111c}
t\mapsto \begin{pmatrix}\dot\sigma_1(t)\\ \dot\sigma_2(t)\end{pmatrix}=\begin{pmatrix}-\eps\omega\sin(\omega t)\\ -\eps\omega\cos(\omega t)\end{pmatrix}
\end{equation}
and we can approximate the waves by piecewise constant functions. These piecewise constant turn out to be of the form \eqref{112} with $\gamma_1=\gamma_2=\eps\omega$. 
Then, plugging the parameters \eqref{ni001} in \eqref{midpointnet} yields, up to $o(\eps^4)$, the value $\delta_m=1.7233\cdot10^{-6}$, whereas $\delta_m^{(N)}=2.0318\cdot10^{-6}$. This amounts to a relative error of the order of $0.15$. Notice that, in the control space, $\delta_m$ corresponds to a square cycle of side $\eps\omega$ (see \eqref{112} with $\gamma_1=\gamma_2=\eps\omega$), whereas $\delta_m^{(N)}$ corresponds to a circular cycle of radius $\eps\omega$ (see \eqref{111c}).
The contribution $\gamma_1\gamma_2\tau^2=(\eps\omega\tau)^2=(\eps\pi)^2/4$ in \eqref{midpointnet} is nothing but the area obtained by integrating the control loop over $[0,4\tau]$. The same contribution in the numerical integration amount to the area of the circle of radius $\eps$ (see \eqref{111b}), which is $\pi\eps^2$: the relative error between these two areas is $0.21$, which is comparable with the relative error $0.15$ between $\delta_m$ and $\delta_m^{(N)}$.


\section{Dependence on the interaction parameter $\lambda$}\label{rem_lambda}
In this section we study briefly the dependence of the equations on the parameter $\lambda=\ln(h/L)/\ln(a/L)$ measuring the strength of the interaction between the two scallops. 
We start, in Subsections~\ref{lambdazero} and~\ref{lambdauno} by studying the limit cases $\lambda=0$ and $\lambda=1$, respectively, to deal with, in Subsection~\ref{lambdarealistic}, with the more realistic case where $\lambda$ is bounded away both from~$0$ and from~$1$.

\subsection{The limit $\lambda=0$}\label{lambdazero}
Assumption \eqref{assumption} implies that $\lambda\in(0,1)$, as already observed, and in the limit as $\lambda\to0^+$ the interaction vanishes, as can be seen both in the expression \eqref{force_density} of the force density and in the expression \eqref{RRR} of the resistance matrix of the system. Indeed (at least formally; the reasoning can be made rigorous by a simple limit process), when $\lambda=0$, and therefore $\Lambda=1$, both $\bef_i^{(j)}(s,t)$ (for $i,j\in\{1,2\}$) are the force densities of a scallop swimming alone in an unbounded fluid, and $\bcR(t;0)$ is a diagonal block matrix, thus 
the equations of motion \eqref{107} are decoupled and correspond to those of two non-interacting scallops.

Notice that the limit $\lambda\to0^+$ can be achieved in two ways. If $h\to L^-$, then $\ln(h/L)\to0^-$ and $\lambda\to0$. This is the case in which the two scallops are sufficiently far away from one another to be considered as non-interacting. This situation violates the second inequality in assumption \eqref{assumption}.
The other case is if $a\to0^+$, which makes the higher order terms in the Resistive Force Theory approximation vanish. Since these are responsible for the interaction \cite{MKL2016}, in this case we would have a system consisting of two non-interacting scallops.

\subsection{The limit $\lambda=1$}\label{lambdauno}
From the definition of $\lambda$, it is immediate to see that $\lambda=1$ if $a=h$ (once more violating assumption \eqref{assumption}). In this case, the distance between the scallops is comparable to their thickness, so that the two swimmers are attached to one another. Because of this and to enforce the non interpenetration of the swimmers, the system is equivalent to one scallop alone: Purcell's Scallop Theorem \cite{Purcell1977} implies that no net displacement can be achieved in this case. 

\subsection{Estimates depending on $\lambda$}\label{lambdarealistic}
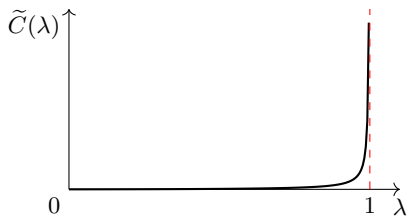
\begin{figure}[!h]
\begin{center}
\begin{tikzpicture}[scale=.8]
\draw[->, color=black] (0,0) -- (5.5,0);
\draw[->, color=black] (0,0) -- (0,3);
\node[text=black,below] at (5,0) {$1$};
\node[text=black,left] at (0,2.75) {$\widetilde C(\lambda)$};
\node[text=black,below] at (5.5,0) {$\lambda$};
\node[text=black,below left] at (0,0) {$0$};
\draw[color=black,thick] plot [samples=200,domain=0:4.985] (\x,{.2*\x*(1+.8*\x)*pow(256*(1-.04*pow(\x,2)),-1)});
\draw[color=red,dashed] (5,0) -- (5,3);
\end{tikzpicture}
\caption{The graph of the function $\lambda\mapsto\widetilde C(\lambda)$.}
\label{Ctilde}
\end{center}
\end{figure}

A non-trivial regime is when the parameter $\lambda$ is far both from~$0$ and from~$1$. To illustrate this case, we provide lower and upper bounds for $\lambda$ in the following relaxed version of assumption \eqref{assumption}: 
\begin{equation}\label{assumption_relaxed}
\kappa a<h<\frac{L}\kappa,
\end{equation}
for a certain $\kappa>0$ to be chosen presently~\footnote{Instead of using $\kappa$ in both sides, one could write \eqref{assumption_relaxed} as $\kappa_1 a<h<L/\kappa_2$, but the relevant estimates would involve the ratio between $\kappa_1$ and $\kappa_2$.}.
Using these inequalities, we obtain that $\lambda\in(\lambda_*(\kappa),\lambda^*(\kappa))$, where
\begin{equation}\label{200}
\!\!\! \lambda^*(\kappa)\coloneqq \frac{-\ln\kappa}{\ln(a/L)}; \; \lambda_*(\kappa)\coloneqq1+\frac{\ln\kappa}{\ln(a/L)}=1-\lambda^*(\kappa).
\end{equation}
The value of $\kappa$ must be chosen in such a way that, \emph{e.g.}, $1/2<\lambda^*(\kappa)<1$ (using the values in \eqref{ni001} for $L$, $a$, $C_\parallel$, and $C_\perp$, the parameter $\kappa$ can be chosen between $2\sqrt{10}$ and $40$).
These bounds on $\lambda$ imply bounds on the constant $C$ in \eqref{C100}, which, in the common approximation $C_\perp= 2C_\parallel$ (see \eqref{ni001} again) are better read in the constant $\widetilde C$ in \eqref{C101}.
Upon noticing that $\lambda\mapsto \widetilde C(\lambda)$ is an increasing function (see Figure~\ref{Ctilde} for a qualitative plot), we obtain that 
\begin{equation}\label{201}
\widetilde C(\lambda_*)<\widetilde C(\lambda)<\widetilde C(\lambda^*);
\end{equation}
the choices of $L$, $a$, $C_\parallel$, and $C_\perp$ in \eqref{ni001} and $\kappa=10$, for instance, yield the values $\widetilde C(\lambda_*)=0.0043$ and $\widetilde C(\lambda^*)=0.0140$.
Estimates of the type \eqref{201} allow us to give an estimate on the displacement \eqref{midpointnet} of the midpoint $\bx_m$: the lower estimate in \eqref{201} provides an estimate on the \emph{minimal} displacement, whereas the maximal displacement yielded by the upper estimate can be overcome by invoking the rate independence of the system, i.e completing the stroke twice as fast makes the swimmer achieve twice the displacement, as expected from the linearity of Stokes flows.

\section{Conclusions}
In this paper, we have studied for the first time the control problem of two scallops hydrodynamically coupled in an unbounded viscous fluid.  The time-reversibility constraint is broken by means of non-local hydrodynamic interactions. We have assumed that the scallops' thickness is non-vanishing, but much smaller than their distance, which, in turn, is much smaller than their length. We have built our control systems on convenient correction of Resistive Force Theory approximation considered in Man et al.~\cite{MKL2016}. To allow analytical progress, we have also assumed the two scallops are approximately parallel to one another in the limit of small angular displacements. This approach has the advantage of working with a finite number of degrees of freedom and sets the basis for future work on discretized interacting micro-swimmers: the position of the hinges between the links and the angle they form are the only parameters used to describe the system, in contrast with the continuum parametrization of a filament \cite{MKL2016,moreau2018asymptotic}.

Our minimal linear model displays the emergence of significant behaviors and regimes that are distinctive of low Reynolds number hydrodynamics. More specifically, we mathematically proved that a particular choice of periodic controls, and therefore of the shape change, provide a non-zero net displacement of the system, overcoming Purcell's scallop theorem.  We have further analyzed how the net displacement varies as a function of the phase difference of two scallops, recovering the maximal displacement predicted previously using other methods (see Theorem~\ref{taylor}). Finally, we have performed an analysis of the displacement as a function of the parameter $\lambda$ associated with the strength of the hydrodynamic interaction between the scallops. Differently from the analysis undertaken in \cite{walker2021control} where the controls used are forces and torques, we choose the velocities of the shape deformations as control functions, with the advantage that only internal actuators of the micro-swimmers are responsible for the motion.
We hope that our analytical solutions will assist and inspire new designs and controls of robotic swimmers that exploits their mutual hydrodynamic interaction in order to propel forwards.

\section*{Acknowledgements}
We are grateful to G.~P.~Alexander and J.~M.~Yeomans for showing us ref.~\cite{AY2008}.
MZ and MM acknowledge that the present research has been partially supported by MIUR (Italian ministry of research) grant Dipartimenti di Eccellenza 2018-2022 (E11G18000350001).

\bibliography{swf}

\end{document}